\documentclass[a4paper,10pt]{article}

\usepackage{amsfonts}
\usepackage{amsmath}
\usepackage{amssymb}

\newtheorem{thm}{Theorem}
\newtheorem{lemma}[thm]{Lemma}
\newtheorem{corol}[thm]{Corollary}
\newtheorem{propos}[thm]{Proposition}
\newtheorem{rmk}[thm]{Remark}

\newtheorem{remark}[thm]{Remark}

\def\E{{\mathbb E}}
\def\P{{\mathbb P}}
\def\R{{\mathbb R}}
\def\Z{{\mathbb Z}}

\def\o{\omega}

\def\O{\Omega}

\def\1{{1\kern-.25em\hbox{\rm I}}}
\def\eu{{1\kern-.25em\hbox{\sm I}}}

\def\R{{\mathbb R}}  
\def\N{{\mathbb N}}  
\def\P{{\mathbb P}}  
\def\Z{{\mathbb Z}}  
\def\C{{\mathbb C}}  
\def\E{{\mathbb E}}  

\let\cal=\mathcal

\def\VV{{\cal V}}

\def\XX{{\cal X}}
\def\X{{\cal X}}

\def\E{{\mathbb E}}

\def\P{{\mathbb P}}
\def\R{{\mathbb R}}
\def\Z{{\mathbb Z}}
\def\N{{\mathbb N}}

\def\C{{\C}}
\def\C{{\cal D}}

\def\sqr{\vcenter{
         \hrule height.1mm
         \hbox{\vrule width.1mm height2.2mm\kern2.18mm\vrule width.1mm}
         \hrule height.1mm}}                  
\def\square{\ifmmode\sqr\else{$\sqr$}\fi}

\def\C{{\bf C}}

\def\be{\begin{equation}}
\def\ee{\end{equation}}

\def\ba{\begin{array}}
\def\ea{\end{array}}
\def\O{{\Omega}}
\def\o{{\omega}}
\def\U{{\Upsilon}}
\def\u{{\upsilon}}

\title{\sf Scaling limit for a drainage network model}
\author{C.~F.~Coletti\thanks{Research supported by FAPESP grant 2006/54511-2}\,\,\footnotemark[4] 
\and 
L.~R.~G.~Fontes
\thanks{Partially supported by CNPq grants 307978/2004-4 and 484351/2006-0, and FAPESP grant 2004/07276-2}\,\,\footnotemark[4] 
\and 
E.~S.~Dias\thanks{Research supported by FAPESP grant 2004/13008-0}\,\,\thanks{Instituto de Matem\'atica
e Estat\'{\i}stica, Universidade de S\~ao Paulo, Rua do Mat\~ao
1010, Ci\-da\-de Universit\'aria, 05508-090 S\~ao Paulo  SP,
Brasil, emails: {\{cristian, lrenato, dias\}@ime.usp.br}}}
\date{}

\begin{document}

\maketitle

\paragraph{Abstract}

We consider the two dimensional version of a drainage network
mo\-del introduced by Gangopadhyay, Roy and Sarkar, and show that
the appropriately rescaled family of its paths converges in
distribution to the Brownian web. We do so by verifying the
convergence criteria proposed by Fontes, Isopi, Newman and
Ravishankar.




\section{Introduction and results}
\renewcommand{\theequation}{1.\arabic{equation}}
\setcounter{equation}{0}
\paragraph{The Two-dimensional Drainage Network Model} Let $\mathbb{Z}$
be the set of the integers and $\O=(\o(z),\, z\in \mathbb{Z}^2)$ be a family of
Bernoulli independent random variables with parameter $p\in(0,1)$.
Let $\mathbb{P}_p$ and $\mathbb{E}_p$ be the probability and expectation induced by those
variables in the product space $(\{0,1\})^{\Z^2}$. Consider  a second family
$\U=(\u(z),\, z\in \Z^2)$ of Bernoulli independent random
variables  with parameter $1/2$. Let $\P_{1/2}$ and $\E_{1/2}$ be the
probability  and expectation induced by these
variables in the product space $(\{0,1\})^{\Z^2}$. Let $\P = \P_p \times
\P_{1/2}$ be the product probability induced by the Bernoulli variables in the
product space $\{0,1\}^{\Z^2} \times \{0,1\}^{\Z^2}$ and let $\E$ be the
expectation operator with respect to this probability.

For $z=(z_1,z_2) \in \mathbb{Z}^2$, say that $z^{\prime}$ is in the next level
if $z^{\prime}=\left(w,z_2+1\right)$ for some $w \in \mathbb{Z}$. Let $h(z)$ be the closest open vertex to
$z$ in the next level with respect to the distance induced by the $l_1$
norm. That it is to say, $\o(h(z))=1$ and $\sum_{i=1}^{2}|h(z)_i -z_i|=\min
\{\sum_{i=1}^{2}|x_i -z_i| : \o(x)=1 \ \mbox{and} \ x_2=z_2 + 1\}$. If there
exist two closest  open vertex to $z$ in the next level, the connection will
be to its left if $\u(z)=0$ and to its right if $\u(z)=1$. Now, let $h^0(z)=z$
and iteratively, for $n \geq 1, h^n(z)=h( h^{n-1}(z))$.

Now let $\mathcal{G}=(V,\mathcal{E})$ be the random directed graph
with vertices $V=\mathbb{Z}^2$ and edges
$\mathcal{E}=\{(u,h(u))\,:\, u \in \mathbb{Z}^2\}$. This model was
proposed in \cite{grs} and will be called the {\em Gangopadhyay, Roy and Sarkar}
(two-dimensional drainage network) {\em model}, or GRS model, for short.

\paragraph{Main result}

It follows from the construction that the GRS model may be seen as
a set of continuous paths, as follows. For any $z=(z_1,z_2)\in
\mathbb{Z}^2$ we define the path $X^z=\{X^z(s),\,s\geq z_2\}$ in
$\R^2$ as the linearly interpolated line composed by all edges
$\{(h^{k}(z),h^{k+1}(z)):k\in\N\}$ of the model, with
$X^z(k)=h^{k}(z)$, $k\in\N$, where $h^k$ is the $k$-th composition
of $h$, $h^0$ meaning the identity. Clearly, $X^z$ is a continuous
path starting at time $z_2$. We Let
\begin{equation}\label{100}\X:=\{X^z:z\in \mathbb{Z}^2\},\end{equation} which we also call
the \emph{drainage network }, and consider its diffusive rescaling
\begin{equation}\label{11}
\X_\delta:=\{(\delta x_1,\delta^2 x_2)\in \R^2:(x_1,x_2)\in
\X\},\end{equation} for $\delta\in(0,1]$. Our main result below
shows that $\X_\delta$ converges in distribution to the Brownian
web.

Several authors constructed random processes that formally
correspond to coalescing one-dimensional Brownian motions starting
from every space-time point \cite{3,4,5,1,1a}. In \cite{1}, the
Brownian web is characterized as a random variable taking values
in a metric space whose points are compact sets of paths. Denote
by $({\cal H},d_{\cal H})$ the complete separable metric space
where the Brownian web is defined. Denote also by ${\cal F}_{\cal
H}$ the corresponding Borel $\sigma$-algebra generated by $d_{\cal
H}$.

The closure in path space of the rescaled drainage network
$\X_\delta$, also denoted by $\X_\delta$, is an $({\cal H},{\cal
F}_{\cal H})$-valued random variable.

\begin{thm} \label{drainageweb}
The rescaled drainage network $\X_\delta$ converges in distribution
to the Brownian web as $\delta \rightarrow 0$.
\end{thm}

In order to prove Theorem~\ref{drainageweb}, we will verify the
convergence criteria of~\cite{1}. To describe them, we need the
following definition. Given $t_0\in\R$, $t>0$, $a<b$, and a
$({\cal H},{\cal F}_{\cal H})$-valued random variable $\VV$, let
$\eta_{_\VV}(t_0,t;a,b)$ be the $\{0,1,2,\ldots,\infty\}$-valued
random variable giving the number of \emph{distinct} points in
$\R\times\{t_0+t\}$ that are touched by paths in $\VV$ which also
touch some point in $[a,b]\times \{t_0\}$.

We can now state the convergence criteria. Let ${\cal D}$ be a
countable dense set of points in $\R^2$.

\begin{thm}[\cite{1}] \label{thmchar}
 {\it Suppose that \( \X_{1}, \X_{2}, \dots \) are \(({\cal H},
{\cal F}_{{\cal H}}) \)-valued random variables with noncrossing
paths. If, in
  addition, the following three conditions are valid, the distribution of
  $\X_n$ converges to the distribution of the standard Brownian
  web.}
\begin{itemize}
\item[($I_1$)] There exist $\theta^y_n\in X_n$ such that for any
    deterministic $y_1,\ldots,y_m\in{\cal D}$,
    $\theta^{y_1}_n,\ldots,\theta^{y_m}_n$ converge in distribution as $n\to
    \infty$ to coalescing Brownian motions (with unit diffusion constant)
    starting at $y_1,\ldots,y_m$.
\item[($B_1$)] $\limsup_{n \to \infty}\sup_{(a,t_0)\in\R^2}
\P(\eta_{_{\X_n}}(t_0,t;a,a+\epsilon)\ge 2)\to0 \hbox{ as }
\epsilon\to0+$; \item[($B_2$)]  $\epsilon^{-1}\limsup_{n \to
\infty}\sup_{(a,t_0)\in\R^2}
\P(\eta_{_{\X_n}}(t_0,t;a,a+\epsilon)\ge 3)\to0\hbox{ as }
\epsilon\to0+$.
\end{itemize}
\end{thm}

The present model can be viewed as a discrete space, long range
version of the two-dimensional Poisson tree introduced
in~\cite{2}. In \cite{2a}, the weak convergence to the Brownian
web of the rescaled collection of paths of that model was
established. The approach is the same as the one pursued here.
There are considerable technical differences in the present work,
though. The main ones are as follows: a) the interaction among the
paths in the present model is a long range one, differently from
the former one; this requires considerably more care in the proofs
of $(I_1)$ and of an estimate on the tail of the distribution of
the coalescence time of two given paths, which we use in the proof
of $(B_2)$ (see Theorems~\ref{Renato1} and~\ref{tailestimate}),
and b) as in \cite{2a}, we resort to an FKG type of argument in an
estimation in the proof of $(B_2)$, and for that we need to
establish a monotonicity property of certain conditional
probabilities on conditioning paths (see
Proposition~\ref{notFKG}); this is also more delicate here than in
\cite{2a}.

Another drainage network for which the rescaled paths were shown
to converge in distribution to the Brownian web is the one
proposed by Scheidegger \cite{sch}, which is equivalent to
coalescing simple random walks. This case was treated in \cite{1}.
In the general context of drainage networks, our results may be
seen as part of the scaling theory for those models. See
\cite{rr}.

The GRS model may be seen also as coalescing random walks, in this
case whose paths are non-crossing and have unbounded dependence
among themselves. The other networks mentioned above may also be
seen as non-crossing coalescing random walks. Coalescing
non-simple random walks which are independent before coalescence,
and whose paths may thus cross one another, were studied in
\cite{nrs} and \cite{bmsv}, and shown under suitable conditions to
converge to the Brownian web when appropriately rescaled.


\section{Proof of Theorem~\ref{drainageweb}}
\renewcommand{\theequation}{2.\arabic{equation}}
\setcounter{equation}{0}

This section is entirely devoted to the proof of Theorem \ref{drainageweb}. Let
us fix a sequence $\delta = \delta_n=1/\sqrt n$ of positive numbers going to zero as $n
\rightarrow \infty$. We want to verify conditions $I_1, B_1$ and $B_2$ of
Theorem \ref{thmchar} for $\X_{\delta_n}$. Due to the translation invariance of the model,
$B_1$ follows from $I_1$ (as will be argued below). So we only have to verify $B_2$ and $I_1$.
This will be done in Subsections~\ref{ssec:b2} and~\ref{ssec:i1} respectively. As a tool for
both this verifications, we start by deriving a bound on the time of coalescence of two paths
of the drainage model.

\paragraph{Coalescing random paths}
Let $u, v \in \mathbb{Z}^2$ be such that $u(1) \leq v(1)$ and
$u(2)=v(2)$. Consider $X^u$ and $X^v$ and, for $t \geq u(2)$, define
\begin{equation}
Z_{t}=Z_{t}(u,v)=X^{v}_{t}-X^{u}_{t}, \label{incremento}
\end{equation}
Notice that $Z_{u(2)} = v(1)-u(1)$.
\begin{rmk}
As argued in~\cite{grs} (see 3.1 in the proof of Theorem 2.1),
$\{Z_{t}(u,v):t\geq u(2)\}$ is a nonnegative martingale in
$\mathbb{L}^{2}$. Also,
$Z_{t}(u,v)\rightarrow 0$ a.s.~as $t\rightarrow\infty$.
\label{martingal}
\end{rmk}

\paragraph{Estimates for the tail of coalescence times} We need to control the tail
of the meeting time of two coalescing random paths starting at the same time at a
distance 1 apart. Define
\begin{equation}
\tau:=\min\{t\geq 0:Z_{t}=0\} \label{eq.tau}
\end{equation}
\begin{thm}\label{tailestimate}
There exists a constant $c_{2}>0$, such that $\mathbb{P}(\tau>t)\leq c_{2}/\sqrt{t}$.
\label{cotatau}
\end{thm}

In order to prove Theorem \ref{tailestimate} it is enough to consider the case
in which $u=(0,0)$ and $v=(1,0)$. By Remark~\ref{martingal}, Skorohod representation holds:
that it is to say, there
exists a Brownian motion with coefficient of diffusion $1$ starting from $1$
and stopping times $0=T_{0},\;T_{1}, ;T_{2},\ldots$ satifying:
\begin{equation}
Z_{t}\stackrel{d}{=}B(T_{t}). \label{skorohold}
\end{equation}
where $0=T_{0},\;T_{1},\;T_{2},\ldots$ are such that
\begin{equation}
T_{t}=\inf\left\{s\geq T_{t-1}:
B(s)-B(T_{t-1})\notin(U_{t}(B(T_{t-1})),V_{t}(B(T_{t-1})))\right\},
\label{eq:tempo}
\end{equation}
where $\{\left(U_{t}(m),V_{t}(m)\right),\;t\geq 1,\;m\in \mathbb{Z} \}$ is a
familiy of random independent vectors and for all
$m\in \mathbb{Z},\;\left(U_{t}(m),V_{t}(m)\right)\stackrel{d}{=}
\left(U_{t^{\prime}}(m),V_{t^{\prime}}(m)\right)$
and  $\left(U_{t}(m),V_{t}(m)\right)\in \{(0,0),A\}$ with
$A=\{\ldots,-2,-1\}\times\{1,2,\ldots\}$.
Note that since $\mathbb{P}\left(Z_{t}\geq 0\right)=1$ we have that
\begin{equation}
\mathbb{P}\left(U_{1}(m)\geq -m\right)=1,\;\forall\;m.
\label{eq.incremento}
\end{equation}
The following result will be needed later on.
\begin{lemma}
For every $p<1$, there exists a constant $c_{1}\in(0,1)$, such that for $m\geq1$
\begin{equation}\label{ub}
\mathbb{P}\left((U_{1}(m),V_{1}(m))=(0,0)\right)\leq c_{1}.
\end{equation}
\label{lem.skorohold}
\end{lemma}

\noindent{\bf Proof}

From the Skorohod representation, one readily concludes that
the left hand side of~(\ref{ub}) can be written as $\P(X^{(m,0)}_1-X^{(0,0)}_1=m)$.
By conditioning on $X^{(0,0)}$, a straightforward computation yields
\begin{equation}\label{ub1}
p^2+\frac{1-q^2}{2(1+q^2)}q^2
\end{equation}
as an upper bound for that probability for all $m$, where $q=1-p$, and this expression is
strictly less than $1$ for $p<1$. $\square$

\medskip

\noindent {\bf{Proof of Theorem \ref{tailestimate}}}

Let $\tau^{\prime}:=\min\{t\geq 0:B(t)=0\}$.
From equations \eqref{skorohold} and \eqref{eq.incremento} we get
\begin{equation}
\mathbb{P}(\tau >t)=\mathbb{P}(\tau^{\prime} > T_{t}). \label{desig.cotatau}
\end{equation}
Now, for $\zeta>0$
\begin{equation}\label{desig.cotatau1}
\mathbb{P}(\tau^{\prime}>T_{t})\leq \mathbb{P}(\tau^{\prime}>\zeta t) + \P(T_{t}\leq\zeta t)
\leq \frac{c}{\sqrt{t}} + \P(T_{t}\leq\zeta t),
\end{equation}
where $c=c(\zeta)\in(0,\infty)$, and the second inequality follows by the well known
behavior of the tail of the distribution of $\tau'$.
By the Markov inequality, we have
\begin{equation}\label{exp.markov}
\mathbb{P}(T_{t}\leq \zeta t)=\mathbb{P}(e^{-\lambda T_{t}}\geq e^{-\lambda\zeta t})
\leq e^{\lambda\zeta t}\,\mathbb{E}\left(e^{-\lambda T_{t}}\right). 
\end{equation}
We now write
\[
T_{t}=\displaystyle\sum_{i=1}^{t}S_{i}(Z_{i-1}),
\]
where $(S_{i}(k),\;i\geq1,\;k\in \mathbb{Z})$ are independent random variables.
(Note that for $k\in \mathbb{Z}$, $\;(S_{i}(k),\;i\geq1)$ are not identically distributed.)
Then, for $\lambda > 0$
\begin{eqnarray}
\mathbb{E}\left(e^{-\lambda T_t}\right) &=&
\mathbb{E}\left[\mathbb{E}\left(\exp\{-\lambda\displaystyle\sum_{i=1}^{t-1}S_{i}(Z_{i-1})\}
\exp\{-\lambda S_{t}(Z_{t-1})\}|\mathcal{F}_{t-1}\right)\right] \nonumber \\
&\leq& \mathbb{E}\left(\exp\{-\lambda\displaystyle\sum_{i=1}^{t-1}S_{i}(Z_{i-1})\}\right)
\displaystyle\sup_{m\in\mathbb{Z}}\mathbb{E}\left(e^{-\lambda S(m)}\right) \nonumber \\
&\leq& \left[\displaystyle\sup_{m\in\mathbb{Z}}\mathbb{E}\left(e^{-\lambda
S(m)}\right)\right]^{t},
\end{eqnarray}
where $\mathcal{F}_t$ is the  $\sigma-$algebra generated by
$\{Z_{0},Z_{1},\ldots,Z_{t}\}$.

Hence,
\begin{equation}
\mathbb{P}(T_{t}\leq\zeta
t)\leq\left[e^{\lambda\zeta}\displaystyle\sup_{m\in\mathbb{Z}}\mathbb{E}\left(e^{-\lambda
S(m)}\right)\right]^{t}, \label{exp.markov1}
\end{equation}
and by Lemma \ref{lem.skorohold} and Skorohod representation, we
get (ommiting $m$'s)
\begin{eqnarray}
\mathbb{E}(e^{-\lambda S})&\leq& \mathbb{P}\left((U,V)=(0,0)\right)\nonumber \\
&+&\mathbb{E} \left(e^{-\lambda S}|(U,V)=(-1,1)\right)
\left(1-\mathbb{P}\left((U,V)=(0,0)\right)\right) \nonumber \\
&=& \left(1-\frac{1}{\cosh(\sqrt{2\lambda})}\right)\mathbb{P}
\left((U,V)=(0,0)\right)+\frac{1}{\cosh(\sqrt{2\lambda})} \nonumber \\
&\leq& \left(1-\frac{1}{\cosh(\sqrt{2\lambda})}\right)c_{1}+\frac{1}{\cosh(\sqrt{2\lambda})}.
\label{exp.markov2}
\end{eqnarray}

Indeed, we first notice that $S$ is the hitting time of $(U,V)$ by
$B$, and thus the case $(U,V)=(-1,1)$ is dominated by the cases
where $(U,V)\ne(0,0)$; this justifies the first inequality
in~\eqref{exp.markov2}. The equality is a well known result (see
e.g.~Theorem 5.7 in~\cite{Du}), and the last inequality follows
from Lemma~\ref{lem.skorohold}. Since $c_1<1$ (uniformly in $m$),
we may find $\lambda_{0}>0$ and $\zeta_{0}>0$ such that
\begin{equation}
c_{3}:=\left(1-\frac{1}{\cosh(\sqrt{2\lambda_{0}})}\right)c_{1}+
\frac{1}{\cosh(\sqrt{2\lambda_{0}})}<e^{-\lambda_{0}\zeta_{0}}.\\
\label{constante}
\end{equation}

Then, from \eqref{exp.markov1}
\begin{equation}
\mathbb{P}\left(T_{t}<\zeta_{0}t\right)\leq
(c_{3}e^{\lambda_{0}\zeta_{0}})^{t}:=c_{4}^{t},
\label{exp.markov0}
\end{equation}
where, from~\eqref{constante}, $c_{4}$ may be taken $<1$. It
follows that there exists $c_{5}\in(0,\infty)$ such that
$c_{4}^{t}\leq c_{5}/\sqrt{t}$. Making now $c_{2}=c+c_{5}$, the
result follows from~\eqref{desig.cotatau1}
and~\eqref{exp.markov0}. $\square$

\subsection{Verification of condition ${\mathbf B_2}$}
\label{ssec:b2}

The drainage network model is translation invariant, then we can
eliminate $\sup_{(a,t_{0})\in\mathbb{Z}^{2}}$ in ($B_2$) and
consider $a=t_{0}=0$.

Then, to prove $(B_{2})$ it suffices to show that
\begin{equation}
\epsilon^{-1}\limsup_{N\rightarrow\infty}\mathbb{P}\left(\eta_{\XX}(0,tN,0,\epsilon\sqrt{N})
\geq3\right)\rightarrow 0, \mbox{\;if\;} \epsilon \rightarrow 0+.
\label{eqfinal}
\end{equation}

For $j\in\mathbb{Z}\cap[0,\epsilon\sqrt{N}]$ let
$X_{j}=\{X^{(j,0)}(k),\,k\in{\mathbb Z}\cap[0,tN] \}$ be the
$\lfloor tN \rfloor$-step trajectory starting in $(j,0)$ .

We now introduce the counting variable $\eta' =
|\left\{X_{j}(tN):1\leq j \leq n \right\}|$, where
$X_{j}(tN)=X^{(j,0)}(\lfloor tN\rfloor)$. Then,
\begin{equation}
\mathbb{P}\left(\eta_{\XX}(0,tN;0,\epsilon\sqrt{N})\geq
3\right)=\mathbb{P}\left(\eta^{\prime}\geq 3\right).
\label{eq.equivalente}
\end{equation}
Since
\[
\left\{\eta^{\prime}\geq3\right\} =  \displaystyle\cup_{j=1}^{n-1}\{X_{j-1}(tN)<X_{j}(tN)<X_{n}(tN)\},
\]
where $n$ is short for $\lfloor \epsilon\sqrt{n}\rfloor$, we have that
\begin{eqnarray} \label{inequalityB2}
&& \mathbb{P}\left(\eta^{\prime}\geq 3 \right)\leq
\displaystyle\sum_{j=1}^{n-1}\mathbb{P}\left(X_{j-1}(tN)<X_{j}(tN)<X_{n}(tN)\right) \nonumber \\
&=&\displaystyle\sum_{j=1}^{n-1}\int_{\bar{\Pi}_{j}}
\mathbb{P}\left(X_{j-1}(tN)<X_{j}(tN)<X_{n}(tN)|X_{j}=\pi\right)\mathbb{P}(X_{j}=\pi)
\nonumber \\
&=&\displaystyle\sum_{j=1}^{n-1}\int_{\bar{\Pi}_{j}}
\mathbb{P}\left(X_{j-1}(tN)<X_{j}(tN)|X_{j}=\pi\right)\nonumber\\
&&\mbox{}\hspace{1cm}\cdot\mathbb{P}\left(X_{j}(tN)<X_{n}(tN)|X_{j}=\pi\right)
\mathbb{P}(X_{j}=\pi),
\label{monotone}
\end{eqnarray}
where $\bar{\Pi}_{j}$ stands for the state space of $X_{j}$. The
last equality folows from the fact that given $X_{j}=\pi$, the
events $\left\{X_{j-1}(tN)<X_{j}(tN)\right\}$ and
$\left\{X_{j}(tN)<X_{n}(tN)\right\}$ are independent.

At this point we want to appeal to the Harris-FKG inequality. So
we must establish monotonicity properties of the conditional
probabilities on the right of~(\ref{monotone}).

We begin by introducing a partial order $\prec$ on $\bar{\Pi}_{j}$
as follows. Given $\pi_1$ and $\pi_2 \in \bar{\Pi}_{j}$, say that
\begin{equation}
\pi_{1}\prec\pi_{2}\Leftrightarrow \pi_{1}(\ell)-\pi_{1}(k)\leq
\pi_{2}(\ell)-\pi_{2}(k), \label{ordem}
\end{equation}
for every $\ell\geq k\geq 0, \, \ell,k \in\mathbb{Z}\cap[0,tN].$

The partial order we have just introduced is an order on the
increments of the trajectories belonging to events in
$\bar{\Pi}_{j}.$ Then, if the increments are increasing, the
corresponding events will also be increasing.

\begin{propos}\label{notFKG}
Let $\pi_1, \pi_2 \in \bar\Pi_j$ be such that $\pi_1 \prec \pi_2$.
Then, for $k<j<n$
\begin{eqnarray} \label{monotone1}
\mathbb{P}\left[X_{k}(tN)<X_{j}(tN)|X_{j}=\pi_1\right] \!\!\!\!\!&\leq&\!\!\!\!\!
\mathbb{P}\left[X_{k}(tN)<X_{j}(tN)|X_{j}=\pi_2 \right],\\
 \label{monotone2}
\mathbb{P}\left[X_{j}(tN) < X_{n}(tN)|X_{j} = \pi_1 \right] \!\!\!\!\!&\geq&\!\!\!\!\!
\mathbb{P}\left[X_{j}(tN) < X_{n}(tN)|X_{j} = \pi_2 \right].
\end{eqnarray}
\end{propos}
Now, using Proposition \ref{notFKG}, and since the increments of
$X_j$ are independent, we may apply Harris-FKG and we find an
upper bound for~(\ref{inequalityB2}) as follows.
\begin{eqnarray}
\nonumber &\displaystyle\sum_{j=1}^{n-1}&\displaystyle\int_{\bar{\Pi}_{j}}
\mathbb{P}\left(X_{j-1}(tN)<X_{j}(tN)|X_{j}=\pi\right)\mathbb{P}(X_{j}=\pi)\\ \nonumber \\
\nonumber &&\cdot\displaystyle\int_{\bar{\Pi}_{j}}
\mathbb{P}\left(X_{j}(tN)<X_{n}(tN)|X_{j}=\pi\right)\mathbb{P}(X_{j}=\pi)\\ \nonumber \\
\nonumber &=&\displaystyle\sum_{j=1}^{n-1}\mathbb{P}\left(X_{j-1}(tN)<X_{j}(tN)\right)
\mathbb{P}\left(X_{j}(tN)<X_{n}(tN)\right)\\ \nonumber \\
\nonumber &\leq&\displaystyle\sum_{j=1}^{n-1}\mathbb{P}\left(X_{j-1}(tN)<X_{j}(tN)\right)
\mathbb{P}\left(X_{0}(tN)<X_{n}(tN)\right)\\ \nonumber \\
\nonumber &\leq&n\mathbb{P}\left(X_{0}(tN)<X_{1}(tN)\right)
\mathbb{P}\left(X_{0}(tN)<X_{\epsilon\sqrt{N}}(tN)\right)\\ \nonumber \\
\nonumber &\leq&
\epsilon\sqrt{N}\mathbb{P}\left(\tau>tN\right)\mathbb{P}\left(\tau_{\epsilon,N}>tN\right),
\end{eqnarray}
where $\tau$ is the coalescing time between the trajectories of
$X^{(0,0)}$ and $X^{(1,0)}$ and $\tau_{\epsilon,N}$ is the
coalescing time for $X^{(0,0)}$ and $X^{(n,0)}$.

From \eqref{eq.equivalente} we get
\begin{equation}
\epsilon^{-1}\mathbb{P}\left(\eta_{\mathcal{\chi}}(0,tN;0,\epsilon\sqrt{N})\geq
3\right)\leq
\sqrt{N}\mathbb{P}\left(\tau>tN\right)\mathbb{P}\left(\tau_{\epsilon,N}>tN\right).
\label{eq.equivalente1}
\end{equation}
From Theorem~\ref{Renato1}, we have that
\[
\limsup_{N\rightarrow
\infty}\mathbb{P}\left(\tau_{\epsilon,N}>tN\right)=\mathbb{P}\left(\tau_{\epsilon,B}>t\right),
\]
where $\tau_{\epsilon,B}$ is the coalescing time between two Brownian motions
starting at a distance
 $\epsilon$ at time zero. Also, recall that
$\mathbb{P}\left(\tau_{\epsilon,B}>t\right)$ is well known to be $O(\epsilon)$. On the other
hand, by Lemma \ref{cotatau}, there exists a constant
$c_{2}^{\prime}>0$ such that
$\mathbb{P}\left(\tau>tN\right)<c_{2}^{\prime}/\sqrt{tN}$. Taking
$\displaystyle\limsup_{N\rightarrow\infty}$ in both sides of
\eqref{eq.equivalente1} we get

\[
\epsilon^{-1}\limsup_{n\rightarrow\infty}\mathbb{P}\left(\eta_{\mathcal{\chi}}(0,tN;0,\epsilon\sqrt{N})\geq
3\right)= O(\epsilon),
\]
which shows \eqref{eqfinal}.

\medskip

\noindent{\bf{Proof of Proposition \ref{notFKG}}}

It is enough to argue~(\ref{monotone1}), since by the symmetry of the model we have
\begin{equation} \label{sym}
\mathbb{P}\left[X_{j}(tN) < X_{n}(tN)|X_{j} = \pi \right]=
\mathbb{P}\left[X_{2j-n}(tN) < X_{j}(tN)|X_{j} = \pi^- \right],
\end{equation}
where $\pi^-$ is the path of $\bar\Pi_j$ whose increments are the
opposite of the respective ones of $\pi$.

The next reduction is that it is enough to consider
$\pi_{1},\pi_{2}$ such that $\pi_{1}(s)=\pi_{2}(s)$ for $0 \leq s<
t_{0}$ and $\pi_{2}(s)=\pi_{1}(s) + 1$ for $s \geq t_{0}$ for some
$t_{0} \in [0,tN]\cap\Z$. From now on, we will be dealing with
such pair of paths $\pi_1$ and $\pi_2$.

\begin{rmk} \label{rmk:condist}
Given that $X_j=\pi$, the distribution of $(\O,\U)$ changes as
follows. Let $I_0=(j,0)$ and for $k\in[1,tN]\cap\Z$, let $I_k$ be
the interval of integer numbers with end points $2\pi(k-1)-\pi(k)$
and $\pi(k)$. ($I_k$ can be a single point.) Then, off
$\cup_{k=0}^{\lfloor tN\rfloor}I_k$, the distribution of $(\O,\U)$
does not change. In $\cup_{k=0}^{\lfloor tN\rfloor}I_k$, the
variables of $(\O,\U)$ are independent of the ones off
$\cup_{k=0}^{\lfloor tN\rfloor}I_k$, and they are also independent
among themselves, except the pairs
$(\o(2\pi(k-1)-\pi(k),k),\u(\pi(k-1),k-1))$ for which $|I_k|>1$,
whose distribution is described below. First, for each
$k\in[1,tN]\cap\Z$, $\o(\pi(k),k)=1$. If $|I_k|>1$, then
$\o(z_1,k)=0$ for $z_1$ in $I_k$ minus its endpoints, and
$\o(2\pi(k-1)-\pi(k),k)$ is a Bernoulli random variable with
parameter $p'=p/(2-p)$; given that $\o(2\pi(k-1)-\pi(k),k)=0$,
then $\u(\pi(k-1),k-1)$ is Bernoulli with parameter $1/2$;
otherwise it is $0$ or $1$ depending on whether
$\pi(k)-\pi(k-1)>0$ or $<0$.

A further point is that given that $X_j=\pi$, the event
$\{X_{k}(tN)<X_{j}(tN)\}$ depends on $\O$ only to the left of and
including $\pi$, and on $\U$ only strictly to the left of $\pi$.
\end{rmk}

In order to prove~(\ref{monotone1}), we couple the
distributions of $(\O,\U)$ given $\pi_1$ and $\pi_2$ conveniently
(in the relevant region of $\Z^2$). Let us denote by
$(\o_i(z),\u_i(z))$, $z\in\Z^2$, the general term of $(\O,\U)$
given $\pi_i$, $i=1,2$, respectively. Then
$(\o_2(z),\u_2(z))=(\o_1(z),\u_1(z))$ for $0\leq z_2<t_0$, and
$(\o_2(z),\u_2(z))=(\o_1((z_1-1,z_2),\u_1(z_1-1,z_2))$ for $0\leq
z_2>t_0$. The coupling on $z_2=t_0$ will be different in different
cases.

\medskip

\noindent{\bf Case 1}: $\pi_1(t_0)<\pi_1(t_0-1)$.

In this case, we make
$(\o_2(z_1,t_0),\u_2(z_1,t_0))=(\o_1(z_1-1,t_0),\u_1(z_1-1,t_0))$
for $z_1\leq\pi_2(t_0)$, and
$(\o_2(z_1,t_0),\u_2(z_1,t_0))=(\o_1(z_1+1,t_0),\u_1(z_1+1,t_0))$
for $z_1>\pi_2(t_0)$.

One readily checks in this case that for every realization
$(\o_1,\u_1)(\cdot)$ of $(\O,\U)$ given $X_j=\pi_1$,
$(\o_2,\u_2)(\cdot)$ given by the above coupling is a realization
of $(\O,\U)$ given $X_j=\pi_2$, and that
\begin{equation}\label{cont}
\{(\o_1,\u_1)(\cdot):\,X_{k}(tN)<X_{j}(tN)\}\subset\{(\o_2,\u_2)(\cdot):\,X_{k}(tN)<X_{j}(tN)\}.
\end{equation}
This then establishes~(\ref{monotone1}) in this case.

\medskip

\noindent{\bf Case 2}: $\pi_1(t_0)\geq\pi_1(t_0-1)$.

We start with a coupling which is the symmetric of the one above.
It does not quite give us~(\ref{cont}), but something nevertheless
useful.

{\bf Auxiliary coupling}:

\noindent Given $(\o_2,\u_2)(\cdot)$, we make
$(\o_1(z_1,t_0),\u_1(z_1,t_0))=(\o_2(z_1+1,t_0),\u_2(z_1+1,t_0))$
for $z_1\geq\pi_1(t_0)$, and
$(\o_1(z_1,t_0),\u_1(z_1,t_0))=(\o_2(z_1-1,t_0),\u_2(z_1-1,t_0))$
for $z_1<\pi_1(t_0)$. One can readily check that this provides a
coupling of $(\O,\U)$ given $X_j=\pi_1$ and $(\O,\U)$ given
$X_j=\pi_2$. It is possible to find realizations of
$(\o_2,\u_2)(\cdot)$ for which with this
coupling~(\ref{monotone1}) does not hold. But the following can be
readily checked.
\begin{equation}\label{cont_aux}
\{(\o_1,\u_1)(\cdot):\,X_{k}(t_0)<X_{j}(t_0)\}
\subset
\{(\o_2,\u_2)(\cdot):\,X_{k}(t_0)<X_{j}(t_0)\}.
\end{equation}
This will be used in the coupling discussed next.

We begin by pointing out that we could not indeed find a direct
coupling for which~(\ref{cont}) holds, so we take an indirect
route. We start by considering $D_0=X_j(t_0)-X_k(t_0)$. In order
to insure the existence of a coupling of $(\o_1,\u_1)(\cdot)$ and
$(\o_2,\u_2)(\cdot)$ satisfying~(\ref{cont}), it is enough to show
that the distribution of $D_0$ given $X_j=\pi_2$ dominates that of
$D_0$ given $X_j=\pi_1$. Indeed, if this holds, then we can find a
coupling of the two distributions such that the domination takes
place almost surely. The respective random variables depend of
course only on the history up to $t_0$, and thus we are free to
choose any coupling of $(\o_1,\u_1)(\cdot)$ and
$(\o_2,\u_2)(\cdot)$ above $t_0$, and choosing as in Case 1, we
readily get~(\ref{cont}).

In order to establish the above mentioned domination, we first
define $\Delta=2(\pi_1(t_0)-\pi_1(t_0-1))$ and
$\tilde D_0=\pi_1(t_0)-\Delta-X_k(t_0)$. Denote
by $D_i$, $i=1,2$, the random variable whose distribution equals
that of the conditional distribution of $\tilde D_0$ given that
$X_j=\pi_i$, $i=1,2$, respectively. Then the above mentioned
domination is equivalent to
\begin{equation}\label{dom}
\P(D_2\geq k-1)\geq\P(D_1\geq k),\,k\geq0.
\end{equation}

We will argue~(\ref{dom}) using different couplings for different
cases. The case $k=0$ of~(\ref{dom}) is established as follows. If
$\pi_1(t_0)=\pi_1(t_0-1)$, there is nothing to prove. If
$\pi_1(t_0)>\pi_1(t_0-1)$, then $D_1\geq0$ is equivalent to the
left hand side of~(\ref{cont_aux}), and $D_2\geq-1$ is equivalent
to $D_2\geq0$, which is in turn equivalent to the right hand side
of~(\ref{cont_aux}). So the auxiliary coupling can be used in this
case. It can also be used if $k=1$. Indeed, in this case,
$D_2\geq0$ is equivalent to the right hand side of~(\ref{cont_aux})
and $D_1\geq1$ is either equivalent to the left hand side
of~(\ref{cont_aux}) (if $\pi_1(t_0)=\pi_1(t_0-1)$), or is in any case
contained in it.

For $k\geq2$, we will resort to another coupling, described as follows.
\begin{enumerate}
    \item For $z_1>\pi_1(t_0)$,
    $(\o_2(z_1,t_0),\u_2(z_1,t_0))=(\o_1(z_1-1,t_0),\u_1(z_1-1,t_0))$;
    \item for  $\pi_1(t_0)-2(\pi_1(t_0)-\pi_1(t_0-1))\leq
    z_1\leq\pi_1(t_0)$, $\o_2(z_1,t_0)=0$ and $\u_2(z_1,t_0)$ are
    i.i.d.~Bernoullis with parameter $1/2$ independent of all else;
    \item for $z_1=\pi_1(t_0)-2(\pi_1(t_0)-\pi_1(t_0-1))-1$,
    $(\o_2(z_1,t_0),\u_2(z_1,t_0))=(\xi\o_1(z_1,t_0),\tilde\u)$,
    where $\xi$ is a Bernoulli with parameter $1/(2-p)$ independent of all else, and $\tilde\u$
    depends only on $\xi\o_1(z_1,t_0)$; if the latter random
    variable vanishes, then $\tilde\u$ is Bernoulli with parameter $1/2$ independent of all
    else; otherwise, $\tilde\u=1$;
    \item and for $z_1<\pi_1(t_0)-2(\pi_1(t_0)-\pi_1(t_0-1))-1$,
    $(\o_2(z_1,t_0),\u_2(z_1,t_0))=(\o_1(z_1,t_0),\u_1(z_1,t_0))$.
\end{enumerate}

One readily checks that for every realization $(\o_1,\u_1)(\cdot)$
of $(\O,\U)$ given $X_j=\pi_1$, $(\o_2,\u_2)(\cdot)$ given by the
above coupling is a realization of $(\O,\U)$ given $X_j=\pi_2$,
and that
\begin{equation}\label{cont1}
\{(\o_1,\u_1)(\cdot):\,D_1\geq
k\}\subset\{(\o_2,\u_2)(\cdot):\,D_2\geq k\}
\end{equation}
in the remainder cases. $\square$

\subsection{Weak convergence to coalescing Brownian motions}
\label{ssec:i1}

The main purpose of this part of the paper is to prove condition
$I1$ of Theorem \ref{thmchar}. Also, we prove condition $B1$,
which is in fact a consequence of condition $I1$ for this model.
\newline
Denote a single trajectory belonging to the drainage network model
and starting from $a\in\mathbb{Z}$ at time $j\in\mathbb{Z}$ by
$X^{a,j} = \{X^{a,j}(l), l \geq j\}$. Then,
\begin{thm}\label{Renato1}
Let $(y_0,s_0),(y_1,s_1),\dots,(y_k,s_k)$ be $k+1$ different points in
$\mathbb{R}^2$ such that $s_0 \leq s_1 \leq \dots \leq s_k$ and if
$s_{i-1}=s_i$ for some $i, i=1,\dots ,k$, then $y_{i-1} < y_i$. If $Z_n^{(i)}
= \{ Z_n^{(i)}(t):= n^{-1/2} X_n^{(i)}(\lfloor nt \rfloor) := n^{-1/2}
X^{\lfloor y_i \sqrt{n}\rfloor,\lfloor s_i n \rfloor}(\lfloor nt \rfloor), t
\geq s_i  \}$; then
\begin{equation}\label{weakconvergence}
\{Z_n^{(i)}, i=0,\dots,k \} \Rightarrow \{W^{(i)}, i=0,\dots,k \},
\end{equation}
as $n \rightarrow \infty$, where the $W^{(i)}$ are $k+1$ coalescing Brownian motions with
positive diffusion coefficient $\sigma$, starting at $\{(y_0,s_0),\dots,(y_k,s_k)\}$.
\end{thm}
\begin{remark}
Here, $\Rightarrow$ stands for weak convergence in $\Pi^{k+1}$ (the product of
$k+1$ copies of $\Pi$, the path space (see~\cite{1}).
\end{remark}
\begin{remark}
$\sigma^2$ is the variance of $X^{0,0}(1)$, which one readily computes from its distribution,
which is in turn straightforward to obtain, as $\frac{q(1+q^2)}{p^2(1+q)^2}$.
\end{remark}

\noindent{\bf Proof of Theorem \ref{Renato1}}

 The proof is divided in three parts. The arguments
are somewhat standard, so we will be sketchy at a few points.

\medskip

{\bf{Part I}}

When $k=0$, the result follows from Donsker' invariance
principle for single paths since $X^{a,j}$ is a random walk with
variance $\sigma^2$.

\medskip

{\bf{Part II}}

In this part, we consider the case where $k \geq 1$ and $s_0 = \dots
= s_k$ (and $y_{i-1} < y_i$ for every $i, i=1, \dots , k$). The general case will
be treated on Part III. Since the
increments are stationary and spatially homogeneous we will take $s_0=y_0=0$.

{\bf{a)}} We first treat the case $k=1$. Without loss of
generality, we suppose that we are working in a probability space
in which
\begin{equation}\label{coup}
Z_n^{(1)} \xrightarrow[n \rightarrow \infty]{a.s.}W^{(1)}.
\end{equation}
The following proposition, besides implying the result in the case
$k=1$, is also the building block of the argument for the other
cases of this part.
\begin{propos}\label{Prop1}
The conditional distribution of $Z_n^{(0)}$ given $Z_n^{(1)}$ converges almost
  surely  to the distribution of $W^{(0)}$ given $W^{(1)}$.
\end{propos}

\noindent{\bf Proof}

The strategy is to approximate the paths $Z_n^{(0)}$ and
$Z_n^{(1)}$ before the time when they first come close (at a
distance of order $n^\alpha$, $\alpha<1/2$) by independent paths;
then to show that after that time they meet quickly.  Let
us write
\begin{equation}\label{111}
X_n^{(i)}(l)=\lfloor y_i\sqrt{n}\rfloor+\sum_{h=1}^{l}S_h^{(i,n)},
\end{equation}
where for each $i=0,1$ and $n\geq1$, $S_h^{(i,n)}, h \geq 1$, are
i.i.d.~random variables distributed as $X^{0,0}(1)$. Next, we
introduce independent copies of $S_h^{(i,n)}, i=0,1,\,n\geq1$. For
each $n\geq1$, let $\{\tilde{S}_h^{(i,n)}, h \geq 1, i=0,1  \}$ be
i.i.d.~random variables such that $\tilde{S}_1^{(0,n)} \sim
X^{0,0}(1)$ and define random variables $\hat{S}_h^{(i,n)}$  by
\begin{equation}\label{112}
\hat{S}_h^{(i,n)} = \left\{\begin{array}{ll}
S_h^{(i,n)}, & \textrm{if $|S_h^{(i,n)}| \leq n^{\alpha}$}\\
\tilde{S}_h^{(i,n)}, & \textrm{otherwise,}\end{array} \right.
\end{equation}
where $\alpha$ is a positive number to be specified later.
Moreover, for $l\geq1$, let $\hat{X}_n^{(i)}(l)=\lfloor
y_i\sqrt{n} \rfloor + \sum_{h=1}^{l} \hat{S}_h^{(i,n)}$. Define
\begin{equation}\label{113}
\hat\tau_n=\min\{l\geq1:\hat{X}_n^{(1)}(l)-\hat{X}_n^{(0)}(l)\leq3n^{\alpha}\}
\end{equation}
and
\begin{equation}\label{114}
\tau_n=\min\{l\geq1:X_n^{(1)}(l)-X_n^{(0)}(l)\leq3n^{\alpha}\},
\end{equation}
and let
\begin{equation}\label{115}
\tilde{X}_n^{(i)}(l)= \left\{\begin{array}{ll}
\hat{X}_n^{(i)}(l), & \textrm{if $l \leq \hat\tau_n$}\\
\hat{X}_n^{(i)}(\hat\tau_n)+\sum_{h>\hat\tau_n}
\check{S}_h^{(i,n)}, & \textrm{otherwise,}\end{array} \right.
\end{equation}
where
\begin{equation}\label{116}
\check{S}_h^{(i,n)}= \left\{\begin{array}{ll}
\tilde{S}_h^{(0,n)}, & \textrm{if $i=0$}\\
S_h^{(1,n)}, & \textrm{if $i=1$.}\end{array} \right.
\end{equation}

From the exponential tail of the distribution of $S_h^{(i,n)}$ (which is
independent of $i$) we readily get the following result.

\begin{lemma}
Given $\beta<\infty$, the complement of the event
\[
A_n := \{ S_h^{(i,n)} = \hat{S}_h^{(i,n)}, i=1,\dots, n^{\beta},i=0,1 \},
\]
has probability going to zero superpolinomially fast as $n\to\infty$.
\end{lemma}
\begin{corol}\label{c}
The event $B_n:= \{ X_n^{(i)}(l) = \hat{X}_n^{(i)}(l); l=1,\dots ,
n^{\beta}; i=0,1 \}$, is such that $\mathbb{P}[B_n^c]\to0$
converges to superpolinomially fast as $n\to\infty$.
\end{corol}

From the construction, the properties stated in the following result are readily verified.
\begin{lemma}
$\tilde{X}_n^{(0)} := \{ \tilde{X}_n^{(0)}(l),\, l \geq 1 \}$ and
$\tilde{X}_n^{(1)} := \{ \tilde{X}_n^{(1)}(l),\, l \geq 1 \}$ are
independent processes and $X_n^{(i)} \sim \tilde X_n^{(i)}$ for
each fixed $i=0, 1$.
\end{lemma}

\begin{corol}\label{c1}
Let $\tilde{Z}_n^{(i)}(t):= n^{-1/2}\tilde{X}_n^{(i)}(\lfloor nt \rfloor),
\, t \geq 0$. Then,
\begin{enumerate}
\item ${\tilde Z}_n^{(1)}\to W^{(1)}$ almost surely as $n
\rightarrow \infty$; \item $(\tilde{Z}_n^{(0)},\tilde{Z}_n^{(1)})
\Rightarrow
  (\tilde{W}^{(0)},W^{(1)})$, where $\tilde{W}^{(0)} \sim
  W^{(0)}$,
  and $\tilde{W}^{(0)}$ and $W^{(1)}$ are independent;
\item the conditional distribution of
$\{\tilde{Z}_n^{(0)}(t),\, t \leq
  \frac{\hat{\tau}_n}{n} \}$ given $\tilde{Z}_n^{(1)}$ converges almost surely
  to that of $\{ \tilde{W}^{(0)}(t),\, t \leq \tau\}$ given
  $W^{(1)}$, where $\tau = \inf \{ t \geq 0 : \tilde{W}^{(0)}(t) = W^{(1)}(t)  \}$.
\end{enumerate}
\end{corol}

\begin{corol}\label{corol1}
The conditional distribution of $\{ Z_n^{(0)}(l), l \leq \tau_n
\}$ given $Z_n^{(1)}$ converges almost surely to that
of $\{ \tilde{W}^{(0)}(t),\, t \leq \tau \}$ given $W^{(1)}$.
\end{corol}

\begin{lemma}\label{l1}
For $\alpha < 1/2$, there exists $\gamma = \gamma (\alpha) < 1$
such that, as $n \rightarrow \infty$,
\begin{equation}
\mathbb{P}[X_n^{(0)}(l) \neq X_n^{(1)}(l),\, l=\tau_n, \dots ,
\tau_n + n^{\gamma}] \rightarrow 0.
\end{equation}
\end{lemma}
\begin{corol}\label{corol2}
$\{Z_n^{(0)}(t),\,t\geq\frac{\tau_n}{n}\}$ converges in probability to
$\{ W^{(1)}(t),\, t \geq \tau \}$.
\end{corol}
Then, Proposition \ref{Prop1} follows directly from Corollaries
\ref{corol1} and \ref{corol2}, and the Markov property.  $\square$

\medskip

\noindent{\bf Proof of Corollary~\ref{c1}}

1. Let $\beta>7/2$ be fixed. We claim that there exists a finite constant $C$ and $\theta>1$ such that
\begin{equation}\label{tail}
 \mathbb{P}[\hat{\tau}_{n} > n^{\beta}]\leq C/n^\theta.
\end{equation}
To argue this we first notice
that $Y^{(n)}:=\hat{X}_n^{(1)}-\hat{X}_n^{(0)}$ is a random walk on the integers starting from
$\lfloor y_i\sqrt{n}\rfloor$ and with increments distributed as
$Y^{(n)}(1)=\hat{S}_1^{(1,n)}-\hat{S}_1^{(0,n)}$. One may readily notice then that
$\hat{\tau}_{n}$ is dominated by the sum of $\lfloor y_i\sqrt{n}\rfloor$ copies of $\tilde\tau_n$,
the hitting time of the negative integers by a random walk with increments distributed as
$Y^{(n)}(1)$, started at the origin. We then have that
\begin{equation}
\P(\hat{\tau}_{n} > n^{\beta})\leq c_1\sqrt{n}\,\P(\tilde\tau_n> c_2\,n^{\beta-1/2})\leq C/n^{(\beta-3/2)/2},
\end{equation}
where the first bound is straightforward from the above domination, and the last one follows from a standard
fact on the tail of the distribution of hitting times for one dimensional random walks
with finite second moment,
like $Y^{(n)}$. The claim is thus justified.

We can then use~(\ref{coup}), Corollary~\ref{c}, and~(\ref{tail}), together with Borel-Cantelli, to conclude.

2. Immediate from Lemma~\ref{l1} and Donsker's theorem.

3. We take a version of $\tilde{Z}_n^{(0)}$, say
  $\bar{Z}_n^{(0)}$, such that
\begin{equation}\label{eq:3}
\left(\bar{Z}_n^{(0)},\tilde{Z}_n^{(1)}\right)
\xrightarrow[\mbox{ }n\rightarrow\infty]{a.s.}
\left(\bar{W}^{(0)},W^{(1)}\right),
\end{equation}
with $(\bar{W}^{(0)},W^{(1)})\sim(\tilde{W}^{(0)},W^{(1)}),$ and
claim that $\bar{\tau}_n \to\bar{\tau}$ almost surely as
$n\to\infty$, where $ \bar{\tau}_n = \min \{ t \geq 0 :
\tilde{Z}_n^{(1)}(t) - \tilde{Z}_n^{(1)}(t) \leq 3 n^{\alpha -
1/2} \} $ and $ \bar{\tau} = \inf \{ t \geq 0 : \bar{W}^{(0)}(t) =
W^{(1)}(t) \}. $ Indeed, the almost sure convergence of
$\bar{\tau}_n$ follows from the property, that with probability
one, two independent Brownian trajectories will instantaneously
cross after touching, i.e.
\[
\mathbb{P}[\, \forall \, \epsilon>0, \bar{W}^{(0)}(t) > W^{(1)}(t)
 \mbox{ for some }  t \leq \bar{\tau} + \epsilon] = 1,
\]
and the following deterministic result (whose proof is an exercise).
The convergence in distribution in the original space follows.
$\square$

\begin{lemma}
Let $f_n, g_n, f, g:\mathbb{R}^+ \to\mathbb{R}$ be continuous
functions such that $f(0) < g(0)$, $T = \inf \{ t \geq 0 :
f(t)=g(t) \}$ is finite and has the property that for every $\delta>0$ 
there exists $t\in[T,T+\delta]$ with $f(t)>g(t)$. Suppose in addition that 
$\lim_{n \rightarrow \infty} \sup_{0 \leq t \leq T + 1} |f_n(t) - f(t)| = 0$, 
and $\lim _{n \rightarrow\infty} \sup _{0 \leq t \leq T + 1} |g_n(t) - g(t)| = 0$, 
and let $T_n = \inf \{ t \geq 0 : g_n(t) - f_n(t) \leq \epsilon_n \}$,
where $(\epsilon_n)$ is a given sequence of numbers 
vanishing as $n\to\infty$. Then, $ \lim_{n \rightarrow \infty} T_n=T$.
\end{lemma}

\noindent{\bf Proof of Corollary~\ref{corol1}}

Immediate from Corollary~\ref{c} and 3.~of Corollary~\ref{c1}.
$\square$

\medskip

\noindent{\bf Proof of Lemma~\ref{l1}}

Let $\sigma_n=\inf\{l\geq\tau_n:\,X_n^{(0)}(l)=X_n^{(1)}(l)\}-\tau_n$.
It can be readily checked that $\sigma_n$ is dominated by $\tilde\sigma_n$,
the hitting time of $0$ by $Z_t((0,0),(\lfloor 3 n^\alpha\rfloor,0))$
(see~(\ref{incremento})). Proceeding as in the proof of Theorem~\ref{cotatau},
we find
\begin{equation}\label{e1}
\P(\tilde\sigma_n>n^\gamma)=\P(\sigma_n'>T_{n^\gamma})
\leq\P(\sigma_n'>\zeta n^\gamma)+\P(T_{n^\gamma}\leq\zeta n^\gamma),
\end{equation}
where $\sigma_n'$ is the hitting time of $0$ by a Brownian motion started at
$\lfloor 3 n^\alpha\rfloor$ (see~(\ref{desig.cotatau}) and~(\ref{desig.cotatau1})).
Choosing now $\gamma$ in $(2\alpha,1)$, one readily checks that the first term on the right
of~(\ref{e1}) vanishes as $n\to\infty$.
That so does the second one follows by~(\ref{exp.markov0}), and the proof is complete.
$\square$

\medskip

\noindent{\bf Proof of Corollary~\ref{corol2}}

By~(\ref{coup}) and Lemma~\ref{l1}, it is enough to show that
for every $\epsilon>0$, as $n\to0$
\begin{equation}\label{c21}
\P(\max_{0\leq l\leq n^\gamma}(X_n^{(1)}(l+\tau_n)-X_n^{(0)}(l+\tau_n))>\epsilon\sqrt{n})\to0.
\end{equation}

By the Markov property and an elementary domination argument, the left hand side of~(\ref{c21})
is bounded above by
\begin{equation}\label{c22}
\P(\max_{0\leq l\leq n^\gamma}\hat Z^{(n)}_l>\epsilon\sqrt{n}),
\end{equation}
where $\hat Z^{(n)}_l=Z_l((0,0),(\lfloor 3 n^\alpha\rfloor,0))$ (see~(\ref{incremento})).
Now by Doob's inequality, (\ref{c22}) is bounded above by constant times
$\E(\hat Z^{(n)}_{n^\gamma})/\sqrt{n}\leq 3 n^{\alpha-1/2}\to0$ as $n\to\infty$. $\square$

\medskip

\noindent{\bf b)}  We now consider the case $k>1$. It suffices to
prove that
\begin{equation} \label{kgeq2}
\mathbb{E}[f_0(Z_n^{(0)})\dots f_k(Z_n^{(k)})] \,
\xrightarrow[\mbox{ }n
      \to\infty
] \, \mathbb{E}[f_0(W^{(0)})\dots f_k(W^{(k)})]
\end{equation}
for any $f_0,\dots,f_k \in C_b(\Pi,\mathbb{R})$, the space of all
real-valued and bounded functions defined on $\Pi$.

The left hand side of~(\ref{kgeq2}) is equivalent to
\begin{equation} \label{kgeq2ii}
\mathbb{E}\{f_1(Z_n^{(1)})
\mathbb{E}[f_0(Z_n^{(0)})f_2(Z_n^{(2)})\dots
f_k(Z_n^{(k)})|Z_n^{(1)}]\}.
\end{equation}
It follows from Lemma~\ref{condindep} below that the conditional
expectation equals
\[
\mathbb{E}[f_0(Z_n^{(0)})|Z_n^{(1)}] \, \,
\mathbb{E}[f_2(Z_n^{(2)})\dots f_k(Z_n^{(k)})|Z_n^{(1)}].
\]
Since $Z_n^{(1)} \xrightarrow{a.s.} W^{(1)}$ by the inductive
hypothesis the conditional expectations above converge, almost
surely, to $\mathbb{E}[f_0(W^{(0)})|W^{(1)}]$ and $\mathbb{E}[f_2
(W^{(2)}) \dots$ $f_k (W^{(k)}) | W^{(1)}],$ respectively. On the
other hand, by the dominated convergence theorem we get the
convergence of~(\ref{kgeq2ii}) to
\begin{equation} \label{kgeq2iii}
\mathbb{E} \{f_1(
  W^{(1)})\mathbb{E}[f_0(W^{(0)})|W^{(1)}] \,
\mathbb{E}[f_2(W^{(2)})\dots f_k(W^{(k)})|W^{(1)}] \}, \nonumber
\end{equation}
which equals
\begin{equation} \label{kgeq2iiii}
\mathbb{E} \{f_1(W^{(1)}) \mathbb{E}[f_0(W^{(0)})f_2(W^{(2)})\dots
f_k(W^{(k)})|W^{(1)}] \}
\end{equation}
Since~(\ref{kgeq2iiii}) equals the right hand side
of~(\ref{kgeq2}), we have the result.  $\square$

\begin{lemma} \label{condindep}
For each $n\geq1$, given $Z_n^{(1)}, Z_n^{(0)}$ is independent of $(Z_n^{(2)},\dots,Z_n^{(k)})$.
\end{lemma}

\noindent{\bf Proof}

$Z_n^{(0)}$ depends only on the variables of $(\O,\U)$ to the left of
and including $Z_n^{(1)}$, and $(Z_n^{(2)},\dots,Z_n^{(k)})$ depends only on the variables of $(\O,\U)$ to the right of
and including $Z_n^{(1)}$. The result follows from the product structure of the conditional distribution of $(\O,\U)$
given $Z_n^{(1)}$
(see Remark~\ref{rmk:condist}), and the facts that, given $Z_n^{(1)}$, the variables of $\O$ along $Z_n^{(1)}$ are fixed,
and when $Z_n^{(0)}$ depends on the variable of $\U$ in a particular position of $Z_n^{(1)}$ (that happens when, and only
when, the increment of $Z_n^{(1)}$ from the given position is strictly positive), then $(Z_n^{(2)},\dots,Z_n^{(k)})$ does
not depend on that particular variable, and vice-versa. $\square$

\medskip

{\bf Part III}

We now take on the general case and proceed by induction on $k\geq0$.
The case $k=0$ was argued in Part I above.
Suppose now that for some $k\geq1$, the result holds up to $k-1$, and
let $(y_0,s_0),\dots,(y_k,s_k)$ be as in the statement of
Theorem~\ref{Renato1}. If the condition on $(s_i)$ of Part II holds, then
there is nothing to prove. Otherwise, there exists $i\in\{1,\dots , k\}$ such
that $s_i > s_{i-1}$. Let $i_0$ be the maximum between such $i$'s.
By the induction hypothesis,
\begin{equation}\label{subvector}
(\{Z_n^{(i)}(t),\,s_i\leq t\leq s_{i_0}\},\,0\leq i<i_0)
\Rightarrow
(\{W^{(i)}(t),\,s_i\leq t\leq s_{i_0}\},\,0\leq i<i_0),
\end{equation}
and we can suppose that such convergence occurs almost surely. By
the Markov property and Part II (strengthened to accomodate the case where
the $y_i$'s may depend on $n$ and converge to, say, $\tilde y_i$ as $n\to\infty$;
this requires minor changes in the argument above), we have that, given~(\ref{subvector}),
\begin{equation} \label{finalconv}
(\{Z_n^{(i)}(t),\,t\geq s_{i_0}\},\,0\leq i\leq k+1)
\Rightarrow
(\{W^{(i)}(t),\,t\geq s_{i_0}\},\,0\leq i\leq k+1).
\end{equation}

Then, Theorem \ref{Renato1} follows from~(\ref{subvector}), (\ref{finalconv}) and
the Markov property. $\square$

\medskip

\noindent{\bf Verification of $B_{1}$}

We note that
$\eta_{\XX_{n}}(0,t;0,\epsilon)$ $\geq 2$ if and only if the two
trajectories starting at time $0$ from the border of the interval
$[0,\epsilon]$ have not met up to time $t$. By Theorem
\ref{Renato1}, the rescaled trajectories of a finite collection of
paths belonging to the drainage network model converge to those of
coalescing Brownian motion. Then, it is straightforward to get,
for the drainage network model, that
\[
\limsup_{n\rightarrow
\infty}\mathbb{P}(\eta_{\XX_{n}}(0,t;0,\epsilon)\geq2)=
2\Phi(\epsilon/\sqrt{2t}\,)-1,
\]
where $\Phi(\cdot)$ is the standard normal distribution function.
This together with the translation invariance of the model
concludes the proof.

\section*{Acknowledgements}
This work started as part of the master's project of the third author and was completed during the
post-doc project of the first one, both at IME-USP and financed by FAPESP. We thank Pablo Ferrari for
discussions on the coupling argument of Subsection~\ref{ssec:b2}.

\end{document}